\documentclass[12pt]{article}
\usepackage{latexsym}
\def\mymedskip{\vskip\medskipamount}
\def\mymedbreak{\par \ifdim\lastskip<\medskipamount
  \removelastskip \penalty-100 \mymedskip \fi}
\def\myaftermedspace{\par \ifdim\lastskip<\medskipamount
  \removelastskip \penalty55\mymedskip\fi}
\newcommand{\eop}{{\unskip\nobreak\hfil\penalty50
          \hskip2em\hbox{}\nobreak\hfil$\Box$
          \parfillskip=0pt \finalhyphendemerits=0 \par}}
\newenvironment{proof}%
{\mymedbreak{\noindent\bf Proof:\enspace}}{\eop\myaftermedspace}
{\mymedbreak{\noindent\bf Proof of Theorem #1:\enspace}}{\eop\myaftermedspace}
\newtheorem{teor}{Theorem}[section]

\newtheorem{lem}[teor]{Lemma}

\newtheorem{prop}[teor]{Proposition}
\newcommand{\beq}{\begin{equation}}
\newcommand{\eeq}{\end{equation}}
\newcommand{\beql}[1]{\begin{equation} \label{#1}}
\newcommand{\eeql}{\end{equation}}
\newcommand{\beqa}{\begin{eqnarray*}}
\newcommand{\eeqa}{\end{eqnarray*}}
\newcommand{\beqal}[1]{\begin{eqnarray} \label{#1}}
\newcommand{\eeqal}{\end{eqnarray}}
\newcommand{\beqan}{\begin{eqnarray}}
\newcommand{\eeqan}{\end{eqnarray}}
\newcommand{\bpf}{\begin{proof}}
\newcommand{\epf}{\end{proof}}

\newcommand{\bB}{{\bf B}}
\newcommand{\bF}{{\bf F}}
\newcommand{\bT}{{\bf T}}

\newcommand{\bb}{\bar{b}}

\newcommand{\Tr}{{\rm Tr}}

\newcommand{\by}{{\bar{y}}}
\newcommand{\ba}{{\bar{a}}}

\begin{document}
\begin{titlepage}
\title{A Class of Permutation Polynomials of $\bF_{2^m}$ Related to Dickson Polynomials }
\date{\today}
\author{Henk D. L.\ Hollmann\\Philips Research Laboratories\\
Prof. Holstlaan 4, 5656 AA Eindhoven\\The Netherlands\\email: {\tt
henk.d.l.hollmann@philips.com}\\
\\Qing Xiang\\Department of Mathematical
Sciences\\University of Delaware\\Newark, DE 19716, USA\\email: {\tt
xiang@math.udel.edu}%
}
\maketitle
\begin{abstract}
We construct a class of permutation polynomials of $\bF_{2^m}$ that are closely related to Dickson polynomials.
\end{abstract}
\end{titlepage}
\newpage

\section{Introduction}

Let $\bF_{q}$ be a finite field of order $q$, where $q$ is a prime
power. We write $\bF_q\setminus\{0\}$ as $\bF_q^*$. A polynomial
$f(X)\in \bF_q[X]$ is called a {\it permutation polynomial} (PP)
of $\bF_q$ if the associated polynomial function $f:c\mapsto f(c)$
from $\bF_q$ to itself is a permutation of $\bF_q$. Permutation
polynomials have been studied extensively in the literature, see
\cite{ln}, \cite{gm1}, \cite{gm2}, \cite{gm3} for surveys of known
results on PPs. A very important class of polynomials whose
permutation behavior is well understood is the class of Dickson
polynomials, which we will define below.

Let $a\in \bF_q$ and let $n$ be a positive integer. We define the {\it
Dickson polynomial} $D_n(X,a)$ over $\bF_q$ by
$$D_n(X,a)=\sum_{j=0}^{\lfloor n/2\rfloor}\frac {n}{n-j}{n-j\choose j}(-a)^jX^{n-2j},$$
where $\lfloor n/2\rfloor$ is the largest integer $\leq n/2$.
Alternatively we may define the Dickson polynomial $D_n(X,a)$ to
be the unique polynomial of degree $n$ over $\bF_q$ such that
\begin{equation}\label{defdickson}
D_n(X+\frac {a}{X},a)=X^n+\left(\frac {a}{X}\right)^n.
\end{equation}
We refer the reader to \cite[p.~8--9]{lmt} for explanations on why
(\ref{defdickson}) can be used to define the Dickson polynomials.
The PPs among the Dickson polynomials have been completely
classified. We state the following theorem due to N\"obauer
\cite{wn}.

\begin{teor}\label{dicksonpp}
The Dickson polynomial $D_n(X,a)$, $a\in \bF_q^*$, is a
permutation polynomial of $\bF_q$ if and only if
$\gcd(n,q^2-1)=1$.
\end{teor}

A proof of this theorem can be found in the original paper of N\"obauer
\cite{wn} or in \cite[p.~356]{ln}. Dickson in his 1896 Ph. D. thesis
observed and partially proved the theorem.

In this note, we construct a family of permutation polynomials of $\bF_{2^m}$. These polynomials are closely
related to Dickson polynomials $D_n(X,1)$ over $\bF_{2^m}$, where $n$ is of the form $2^k-1$. (See
Proposition~\ref{Hdickson} for the relation.) We state our main results as follows.

Let $m\geq1$ be an integer, let $k$ be an integer in $\{1,\ldots,
m-1\}$ with $\gcd(k,m)=1$, and let $r\in\{1,\ldots, m-1\}$ be such
that $kr\equiv 1 \bmod m$. Define the integer $m'$ by $kr=1+m m'$
and write $q=2^m$ and $\sigma=2^k$. Throughout the rest of the
note, we will keep the definitions of $m,k,r,m',q,\sigma$ fixed.
We will use $\Tr$ to denote the trace from $\bF_q$ to $\bF_2$ and
for $e\in\bF_2$ we set
$$\bT_e=\{x\in \bF_q\mid \Tr(x)=e\}.$$
Also we define $\Tr(X)$ to be the following polynomial in $\bF_2[X]$.
$$\Tr(X):=X+X^2+\cdots+X^{2^{m-1}}.$$

For $\alpha, \gamma$ in $\{0,1\}$, we define the polynomials
\[H_{\alpha,\gamma}(X):= \gamma \Tr(X) + \frac{\left(\alpha \Tr(X) + \sum_{i=0}^{r-1} X^{\sigma^i}\right)^{\sigma+1}} {X^2}.\]
(Note that $H_{\alpha,\gamma}(X)$ is indeed a polynomial in $X$
with coefficients in $\bF_2$ and $H_{\alpha,\gamma}(0)=0$.)

Our main theorem is

\begin{teor}\label{mainthm} Let $m, k$ be positive integers with $\gcd(k,m)=1$,
let $r\in\{1,\ldots, m-1\}$ be such that $kr\equiv 1 \;(\bmod
\;m)$, and let $\alpha, \gamma\in \{0,1\}$. Then the mapping $H_{\alpha, \gamma}: x\mapsto H_{\alpha,\gamma}(x)$, $x\in\bF_q$, maps $\bT_0$ bijectively to $\bT_0$, and maps $\bT_1$ bijectively to $\bT_{r+(\alpha+\gamma)m}$. In particular, the polynomial $H_{\alpha,\gamma}(X)$ is a PP of $\bF_{2^m}$ if and only if $r+(\alpha+\gamma)m \equiv 1$ {\em (mod 2)}.
\end{teor}

The polynomials $H_{\alpha,\gamma}(X)$ arose in our recent work on the association scheme afforded by the action
of ${\rm PGL}(2,q)$ on the set of exterior lines to a non-degenerate conic in ${\rm PG}(2,2^m)$ \cite{hxpc}. In
order to prove that the fusion by the Frobenius map of the aforementioned association scheme is pseudocyclic, we
need to investigate the permutation behavior of the polynomials $H_{\alpha,\gamma}(X)$. We believe that the
polynomials $H_{\alpha,\gamma}(X)$ are of independent interest. In Section 2, we will explain the connection
between Dickson polynomials and the polynomials $H_{\alpha,\gamma}(X)$. In Section 3, we give a proof of our main
theorem.

\section{Relating $H_{\alpha,\gamma}(X)$ to Dickson Polynomials}

Let $m,k,r,m',q,\sigma$ be defined as in Section 1, so that
$\gcd(k,m)=1$ and $kr=1+m'm$. For $\alpha, \beta\in\{0,1\}$, we
define the polynomials
\[ f_\alpha(X):= \alpha \Tr(X) + \sum_{i=0}^{r-1} X^{\sigma^i}, \]
and
\[ g_\beta(X):= \beta \Tr(X)+\sum_{j=0}^{k-1} X^{2^j}. \]
We will use $f_{\alpha}$ and $g_{\beta}$ to denote the associated
polynomial functions from $\bF_q$ to $\bF_q$. Also using
$f_{\alpha}(X)$, we can rewrite $H_{\alpha,\gamma}(X)$ as
$$H_{\alpha,\gamma}(X)=\gamma \Tr(X)+\frac {f_{\alpha}(X)^{\sigma
+1}} {X^2}.$$ In the following lemma we collect the properties of
$f_\alpha$ and $g_\beta$ that will be used in the sequel. Most of
the properties are straightforward and appeared in \cite{dob}. For
completeness, we provide a (different) proof here.

\begin{lem}\label{fgprop}
The maps $f_\alpha$ and $g_\beta$ are both linear on $\bF_q$. They
also have the following additional properties.\\
{\em (i)} For every $x\in\bF_q$, we have $\Tr(f_\alpha(x))=(r+\alpha m) \Tr(x)$ and $f_\alpha(1)=r+\alpha m$.\\
{\em (ii)} For every $x\in\bF_q$, we have $\Tr(g_\beta(x))=(k+\beta m) \Tr(x)$ and $g_\beta(1)=k+\beta m$.\\
{\em (iii)} For every $x\in\bF_q$, we have $f_\alpha(x)^{\sigma}+f_\alpha(x)= x^2+x$ and $g_\beta(x)^2+g_\beta(x)=x^\sigma+x$.\\
{\em (iv)} We have that $f_\alpha$ maps $\bT_0$ bijectively onto
$\bT_0$ and maps $\bT_1$ bijectively onto $\bT_{r+\alpha m}$. In
particular, $f_\alpha$ is a permutation on $\bF_q$ if and only if
$r+\alpha m \equiv 1$ {\em (mod 2)}.\\
{\em (v)}
We have that $g_\beta$ maps $\bT_0$ bijectively onto $\bT_0$ and maps $\bT_1$ bijectively
onto $\bT_{k+\beta m}$. In particular, $g_\beta$ is a permutation on $\bF_q$ if and only if
$k+\beta m \equiv 1$ {\em (mod 2)}.\\
{\em (vi)} For every $x\in\bF_q$, we have
$f_\alpha(g_\beta(x))=g_\beta(f_\alpha(x))= x+\delta \Tr(x)$ with
\[\delta = m'+\alpha k +\beta r +\alpha \beta m.\]
We have that $1+\delta m = (r+\alpha m)(k+\beta m)$.\\
{\em (vii)} For every $y\in\bF_q$ and for every $\lambda\in\bF_2$, we have
$g_\beta(y)=g_0(\by) + \theta \Tr(y)$ with $\by=y+\lambda \Tr(y)$ and $\theta = \beta + \lambda k$.
Here, the element $\theta$ of $\bF_2$ satisfies $m\theta = k+\beta m +k(1+\delta m)$.
\end{lem}
\begin{proof}
The claims (i) and (ii) are trivial. (Simply note that
$\Tr(1)=m$.) The claims in (iii) are easily verified.

Since by (i) $f_\alpha$ is linear, maps $\bT_0$ to $\bT_0$, and maps $\bT_1$ to $\bT_{r+\alpha m}$, claim (iv) is
equivalent to the claim that if $f_\alpha(x)=0$ and $\Tr(x)=0$, then $x=0$. To show this, suppose that
$f_\alpha(x)=0$. By (iii), we have that $0=f_\alpha^\sigma(x) + f_\alpha(x)=x^2+x$, so $x=0$ or $x=1$. Now
$\Tr(1)=m$ and by (i) we have that $f_\alpha(1)=r+\alpha m$. So $\Tr(1)=0$ and $f_\alpha(1)=0$ would imply that
$r\equiv m\equiv 0\bmod 2$, contradicting the assumption that $rk\equiv 1 \bmod m$.

Similarly, claim (v) is equivalent to the claim that if $g_\beta(x)=0$ and $\Tr(x)=0$, then $x=0$, which can be
shown in the same way as the claim for $f_\alpha$ above. Indeed, suppose that $g_\beta(x)=0$. Then by (iii) we
have that $0=g_\beta^2(x)+g_\beta(x)=x^\sigma+x$, and since $\sigma=2^k$ with $\gcd(k,m)=1$, we conclude that
$x=0$ or $x=1$. Again the assumptions that $m=\Tr(1)\equiv 0 \bmod 2$ and $k+\beta m = g_\beta(1) \equiv 0\bmod
2$ would imply that $m\equiv k\equiv 0\bmod2$, which contradicts $rk\equiv 1 \bmod m$.

To prove claim (vi), first note that
\beqa
f_0(g_0(x)) = g_0(f_0(x))
    &=& \sum_{i=0}^{r-1}\sum_{j=0}^{k-1} x^{2^{ki+j}}\\
    &=&\sum_{t=0}^{kr-1} x^{2^t}\\
    &=& x^{2^{m'm}} + \sum_{t=0}^{m'm-1} x^{2^t}\\
    &=& x + m' \Tr(x).
\eeqa
Then use (i), (ii), and the definitions of $f_\alpha$ and $g_\beta$
to verify claim (vi) for arbitrary $\alpha$ and $\beta$.
The last part of claim (vi) follows immediately from (i) and (ii) by taking $x=1$.

Finally, let $\by=y+\lambda\Tr(y)$. Then
\beqa g_\beta(y) &=& \beta \Tr(y) + g_0(y)\\
        &=& \beta \Tr(y) +g_0(\by+\lambda\Tr(y)) \\
        &=&  \beta \Tr(y) +k\lambda\Tr(y) +g_0(\by)\\
        &=& g_0(\by) +\theta \Tr(y)
\eeqa
with $\theta = \beta + k\lambda$. To prove the last part of (vii), simply take $y=1$.
This completes the proof of the lemma.
\end{proof}

In what follows, we will show that the polynomial
$H_{\alpha,0}(X)$ is closely related to Dickson polynomials. First
we observe that in characteristic 2, the Dickson polynomials
$D_{2^k-1}(X,1)$ over $\bF_q$ are closely related to the
linearized polynomial

$$T_k(X)=X+X^2+\cdots +X^{2^{k-2}}+X^{2^{k-1}}$$

To simplify notation, we will use $D_n(X)$ to denote $D_n(X,1)$ over $\bF_{q}$.

\begin{prop}\label{dicksonandlin}
For any $k\geq 1$, $D_{2^k-1}(X)=X^{2^k+1}T_k(1/X)^2.$
\end{prop}

This proposition can be proved by induction, see \cite{cm}.

We are now ready to relate
$H_{\alpha,0}(X)=f_\alpha(X)^{\sigma+1}/X^2$ to $D_{2^k-1}(X)$. We
state our result in the following proposition.

\begin{prop}\label{Hdickson}
Let $m,k,r,m',q,\sigma$ be given as in Section 1 with
$\gcd(k,m)=1$. Let $\alpha\in \{0,1\}$ be such that $r+\alpha
m\equiv 1$ {\em (mod 2)}, and let $\beta\in\{0,1\}$ be defined by
$\beta\equiv m'+\alpha k$ {\em (mod 2)}. Then for every $x\in\bF_q^*$, we
have
$$H_{\alpha,0}(g_{\beta}(x))=x^{\sigma +1}/g_{\beta}(x)^2= \left\{ \begin{array}{ll}
    1/D_{2^k-1}(1/x) & \mbox{if $\beta=0$,} \\
    1/D_{2^{m-k}-1}(1/x)^{2^k}& \mbox{if $\beta=1$.}
            \end{array}
        \right. $$
In particular, $H_{\alpha,0}(X)$ is a PP of $\bF_q$ if and only if
$r+\alpha m\equiv 1$ {\em (mod 2)}.
\end{prop}

\bpf First note that if $m\equiv 0$ (mod 2) then $k\equiv 1$ (mod 2) and $r\equiv 1$ (mod 2). So it is always
possible to choose $\alpha$ such that $r+\alpha m\equiv 1$ (mod 2). When $r+\alpha m\equiv 1$ (mod 2), then by
Lemma~\ref{fgprop}, part (iv), the linear map $f_{\alpha}$ is a permutation of $\bF_q$. Its inverse is
$$g_{\beta}(x)=x+x^2+\cdots +x^{2^{k-1}}+\beta \Tr(x)=T_k(x)+\beta \Tr(x),$$
where $\beta$ is defined in the statement of the proposition. In particular, we have $k+\beta m\equiv 1$ (mod 2),
by Lemma~\ref{fgprop}, part (v). Therefore for $x\in\bF_q^*$ we have \beq\label{Eqs}
H_{\alpha,0}(g_{\beta}(x))=f_\alpha(g_\beta(x))^{\sigma +1}/g_\beta(x)^2= x^{\sigma +1}/g_{\beta}(x)^2. \eeq

\noindent{\bf Case 1.} $\beta=0$. (Hence $k$ is odd.) In this
case, $g_{\beta}(x)=T_k(x)$. Therefore, for every $x\in\bF_q^*$,
by Proposition~\ref{dicksonandlin} and (\ref{Eqs}), we have
\begin{equation}\label{h0dickson}
H_{\alpha,0}(g_{\beta}(x))=x^{\sigma
+1}/T_k(x)^2=1/D_{2^k-1}(1/x).
\end{equation}
Since $\gcd(k,m)=1$ and $k$ is odd, we see that $\gcd(2^k-1, q^2-1)=2^{\gcd(k, 2m)}-1=1$. By
Theorem~\ref{dicksonpp}, $D_{2^k-1}(X)$ is a PP of $\bF_q$. So (\ref{h0dickson}) implies that $H_{\alpha,0}(X)$
is a PP of $\bF_q$.

\noindent{\bf Case 2.} $\beta=1$. (Hence $k+m\equiv 1$ (mod 2).)
In this case, for $x\in \bF_q^*$ (so $g_\beta(x)\neq 0$), we
have

\begin{eqnarray*}
H_{\alpha,0}(g_{\beta}(x))&=&\frac {x^{2^k+1}} {x^2+x^{2^2}+\cdots +x^{2^k}+(x+x^2+\cdots +x^{2^{m-1}})} \\
      &=&\frac {x^{2^k+1}} {x^{2^k+1}+\cdots +x^{2^{m}}} \\
      &=&\left(\frac {x^{2^{m-k}+1}} {x^2+x^{2^2}+\cdots +x^{2^{m-k}}} \right)^{2^k}\\
      &=&\left(1/D_{2^{m-k}-1}(1/x)\right)^{2^k}\\
\end{eqnarray*}
Since $\gcd(m-k,m)=1$ and $m-k$ is odd, we see that $\gcd(2^{m-k}-1,q^2-1)=1$.
Hence by Theorem~\ref{dicksonpp}, $H_{\alpha,0}(X)$ is a PP of $\bF_q$.

Finally if $H_{\alpha,0}(X)$ is a PP of $\bF_q$, then $r+\alpha
m\equiv 1$ (mod 2) since $H_{\alpha,0}(0)=0$ and
$H_{\alpha,0}(1)=r+\alpha m$.

This completes the proof.
\epf

\section{Proof of the Main Theorem}

We will need the following lemmas in the proof of our main theorem.

\begin{lem}\label{Hprop} With the definitions of $H_{\alpha,\gamma}(X)$
and $f_{\alpha}(X)$ given in Section 1 and 2, for $x\in\bF_q^*$, we
have
\[H_{\alpha,\gamma}(x)=
\gamma \Tr(x) + \biggl(\frac {f_\alpha(x)} {x}\biggr)^2+\frac{f_\alpha(x)}{x}+f_{\alpha}(x),\]
and
\[\Tr(H_{\alpha,\gamma}(x))=\left(r+(\alpha+\gamma)m\right) \Tr(x).\]
\end{lem}

\begin{proof}
The first assertion follows from  part (iii) of Lemma~\ref{fgprop}.
Indeed, for $x\in\bF_q^*$, since $f_{\alpha}(x)^\sigma=f_\alpha(x)+x^2+x$,
we have that $f_{\alpha}(x)^{\sigma+1}=f_\alpha(x)^2+f_\alpha(x)(x^2+x)$. Now
\begin{eqnarray*}
H_{\alpha,\gamma}(x)&=&\gamma \Tr(x)+\frac {f_{\alpha}(x)^{\sigma +1}} {x^2}\\
&=&\gamma \Tr(x)+\biggl(\frac {f_\alpha(x)}
{x}\biggr)^2+\frac{f_\alpha(x)}{x}+f_{\alpha}(x).
\end{eqnarray*}

The second assertion follows from the first one in combination
with part (i) of Lemma~\ref{fgprop}.
\end{proof}

Let $m,k,q,\sigma$ be as before. Define $\bB_0=(\bF_q\setminus\{1\})\cup \{\infty\}$ and $\bB_1=\{z\in
\bF_{q^2}\setminus \{0,1\} \mid z^q=z^{-1}\}$. Note that $\bB_1=\{\theta^{(q-1)i} \mid i=1, \ldots, q\}$, for a
primitive element $\theta$ of $\bF_{q^2}$. Also, define the map $\phi$ from $\bF_{q^2}\cup\{\infty\}$ to itself
by
\beq\label{phidef} \phi(z)=1/(z+z^{-1}), \eeq
where the usual convention on the symbol $\infty$ is adopted (in
particular, $\phi(0)=\phi(\infty)=0$ and $\phi(1)=\infty$). Finally, define the maps $w_0$ and $w_1$ from
$\bF_{q^2}\cup\{\infty\}$ to itself by \beq\label{wdef} w_0(z)=z^{\sigma-1}, \qquad w_1(z)=z^{\sigma+1}, \eeq for
$z\in \bF_{q^2}$ and, in addition, $w_e(\infty) = \infty$ for $e=0,1$. Our interest in the sets $\bB_e$ and the
maps $\phi$, $w_0$, and $w_1$ is explained by the following lemma (see also \cite{note}, Lemma~1).

\begin{lem}\label{perm}
{\em (i)} For $e\in\bF_2$, the map $\phi$ maps $\bB_e$ two-to-one onto
$\bT_e$.\\
{\em (ii)} The map $w_0$ is a permutation of $\bB_0$, and it is a permutation of $\bB_1$ if and only if $k\equiv 1$ (mod 2).\\
{\em (iii)} The map $w_1$ is a permutation of $\bB_0$ if and only if $m\equiv 1$ (mod 2), and a permutation of
$\bB_1$ if and only if $m+k\equiv 1$ (mod 2).
\end{lem}
\bpf (i). Let $u: \bF_{q^2}\cup\{\infty\}\rightarrow \bF_{q^2}\cup\{\infty\}$ be defined by $u(z)=1/(z+1)$. Then
the map $u$ is one-to-one from $\bF_{q^2}\cup\{\infty\}$ to itself, and $\phi(z)=u(z)^2+u(z)$ for all
$z\in\bF_{q^2}\cup\{\infty\}$. If $z\in \bB_0$, then $u(z)\in\bF_q$, hence $\phi(z)\in\bT_0$. Since $u$ maps
$\bB_0$ bijectively to $\bF_q$, and the map $z\mapsto z^2+z$ is two-to-one from $\bF_q$ to $\bT_0$, we see that
$\phi$ is two-to-one from $\bB_0$ to $\bT_0$. Now if $z\in\bB_1$, then $\phi(z)^q=\phi(z)$, hence
$\phi(z)\in\bF_q$. But $u(z)\notin \bF_q$, so $\phi(z)\in\bF_q\setminus\bT_0 =\bT_1$. One can further verify that
$u$ maps $\bB_1$ bijectively to the set $\{x\in \bF_{q^2}\mid x^q=x+1\}$, which is mapped two-to-one onto $\bT_1$
by the map $z\mapsto z^2+z$. This shows that $\phi$ is two-to-one from $\bB_1$ to $\bT_1$.

(ii) and (iii). For any integer $s$, the map $z\mapsto z^s$ is a permutation of $\bB_0$ if and only if
$\gcd(s,2^m-1)=1$, and a permutation of $\bB_1$ if and only if $\gcd(s,2^m+1)=1$. Now suppose that $\gcd(k,m)=1$.
If $s=2^k-1$, then $\gcd(s,2^m-1)=1$ (hence $z\mapsto z^{2^k-1}$ is a permutation of $\bB_0$),  and
$\gcd(s,2^m+1)=\gcd(2^k-1, 2^{2m}-1)/\gcd(2^k-1,2^m-1)=2^{\gcd(k,2m)}-1$. So $\gcd(2^k-1, 2^{m}+1)=1$ if and only
if $k$ is odd. Hence the map $z\mapsto z^{2^k-1}$ is a permutation of $\bB_1$ if and only if $k\equiv 1 \bmod 2$.
Next, if $s=2^k+1$, then $\gcd(s,2^m-1) =\gcd(2^k+1,2^m-1)=\gcd(2^{2k}-1, 2^m-1)/\gcd(2^k-1,
2^m-1)=2^{\gcd(m,2k)}-1$. So $\gcd(2^k+1, 2^m-1)=1$ if and only if $m$ is odd. Finally
\beqa \gcd(2^m+1, s) &=& \gcd(2^m+1, 2^k+1)\\
        &=& \gcd(2^{2m}-1,2^k+1)/\gcd(2^m-1,2^k+1)\\
        &=&(2^{\gcd(2m,2k)}-1)(2^{\gcd(m,k)}-1)/((2^{\gcd(m,2k)}-1)(2^{\gcd(2m,k)}-1)),
\eeqa
so $\gcd(2^m+1, 2^k+1)=1$ if and only if precisely one of $k$, $m$ is odd.
\epf

In the sequel we will use the map $\phi$ defined in (\ref{phidef})
and the maps $w_0$ and $w_1$ defined in (\ref{wdef}) to simplify
an equation involving $g_{\beta}(x)$, $x\in \bF_{q^2}$, using the
following lemma.

\begin{lem}\label{zsumexp}
Let $m,k,q,\sigma$ be defined as in Section 1.\\
{\em (i)} For $z\in\bF_{q^2}\setminus\{0,1\}$,
we have that
\[ \sum_{j=1}^{k}(z+z^{-1})^{-2^j} = (z^{\sigma-1}+z^{1-\sigma})/(z+z^{-1})^{\sigma+1}.\]
{\em (ii)} For $z\in\bF_{q^2}\setminus\{0,1\}$, we have that
$g^2_0(\phi(z))=(\phi(z))^{\sigma+1} / \phi(w_0(z))$ and
$1+g^2_0(\phi(z)) = (\phi(z))^{\sigma+1} / \phi(w_1(z))$.
\end{lem}
\begin{proof}
To prove (i), we use induction on $k$. For $k=1$, we have
$\sigma=2$ and the assertion is trivial. Next, if the assertion
holds for $k$, then using induction hypothesis, we have for
$z\in\bF_{q^2}\setminus\{0,1\}$
\beqa \sum_{j=1}^{k+1}
(z+z^{-1})^{-2^j} &=&
(z^{\sigma-1}+z^{1-\sigma})/(z+z^{-1})^{\sigma+1}
            + (z+z^{-1})^{-2\sigma} \\
    &=& \biggl((z^{\sigma-1}+z^{1-\sigma})(z+z^{-1})^\sigma+(z+z^{-1})\biggr)
            /(z+z^{-1})^{2\sigma+1} \\
    &=& (z^{2\sigma-1}+z^{1-2\sigma})/(z+z^{-1})^{2\sigma+1},
\eeqa and the assertion holds also for $k+1$. This proves (i).

The first assertion in (ii) is a direct consequence of (i); the second assertion is a consequence of the fact that $(z+z^{-1})^{\sigma+1} =
z^{\sigma+1} + z^{\sigma-1}+z^{-\sigma+1} +z^{-\sigma-1}$.
\end{proof}

We are now ready to give the proof of our main theorem.
\vspace{0.1in}

\noindent{\bf Proof of Theorem~\ref{mainthm}.} By
Lemma~\ref{Hprop}, the mapping $H_{\alpha, \gamma}: x\mapsto
H_{\alpha,\gamma}(x)$, $x\in\bF_q$, maps $\bT_0$ to $\bT_0$, and
maps $\bT_1$ to $\bT_{r+(\alpha+\gamma)m}$. So it suffices to show
that $H_{\alpha,\gamma}$ is injective on both $\bT_0$ and $\bT_1$.
For $x\in\bT_0$, we have $H_{\alpha, \gamma}(x)=H_{0, 0}(x)=H_{1,
0}(x)$. Since $\gcd(m,r)=1$, it is always possible to choose
$\alpha\in\{0,1\}$ such that $r+\alpha m\equiv 1$ (mod 2). It
follows from Proposition~\ref{Hdickson} that $H_{\alpha, \gamma}$
maps $\bT_0$ to $\bT_0$ bijectively.

Now we show that if
\begin{equation}\label{1-to-1equ1}
H_{\alpha, \gamma}(x)=H_{\alpha, \gamma}(y), \;{\rm and}\; x,y\in\bT_1,
\end{equation}
then $x=y$.
Simplifying (\ref{1-to-1equ1}), we get
\begin{equation}\label{1-to-1equ2}
\frac{f_{\alpha}(x)^{\sigma +1}}{x^2}=\frac{f_{\alpha}(y)^{\sigma +1}}{y^2}.
\end{equation}
Since $\gcd(m,k)=1$, it is possible to choose $\beta\in\{0,1\}$
such that $k+\beta m\equiv 1$ (mod 2). By Lemma~\ref{fgprop}, part
(v), $g_{\beta}$ maps $\bT_1$ bijectively to $\bT_{k+\beta
m}=\bT_1$. Let $a,b$ be elements of $\bT_1$ such that $g_{\beta}(a)=x$ and
$g_{\beta}(b)=y$. Substituting $x,y$ in (\ref{1-to-1equ2}) by
$g_{\beta}(a)$ and $g_{\beta}(b)$ respectively, and applying
Lemma~\ref{fgprop}, part (vi), we have
\begin{equation}\label{1-to1equ3}
\frac {(a+\delta)^{\sigma +1}} {g_{\beta}(a)^2}=\frac {(b+\delta)^{\sigma +1}} {g_{\beta}(b)^2},
\end{equation}
where $\delta\equiv m'+\alpha +\beta r$ (mod 2). Set
$a+\delta=\ba$ and $b+\delta=\bb$. Applying Lemma~\ref{fgprop},
part (vii), with $\lambda=\delta$, we have
\begin{equation}\label{1-to1equ4}
\frac {\ba^{\sigma +1}} {g_{0}(\ba)^2+\theta}=\frac {\bb^{\sigma
+1}} {g_{0}(\bb)^2+\theta},
\end{equation}
where $\theta\equiv \beta+\delta k$ (mod 2). Note that since
$k+\beta m\equiv 1$ (mod 2), we have
\begin{equation}\label{congr}
m\theta\equiv 1+k(1+\delta m) \; ({\rm mod} \;2)
\end{equation}
If one of $\ba, \bb$ is zero, then by (\ref{1-to1equ4}), the other
one is also zero, hence $a=b$, and therefore $x=y$. So from now on
we assume that $\ba\neq 0$ and $\bb\neq 0$. Then we obtain from
(\ref{1-to1equ4}) that
\begin{equation}\label{1-to1equ5}
\frac {g_{0}(\ba)^2} {\ba^{\sigma +1}}+ \frac {\theta}
{\ba^{\sigma +1}}=\frac {g_{0}(\bb)^2} {\bb^{\sigma +1}}+\frac
{\theta} {\bb^{\sigma +1}}.
\end{equation}

Note that $\Tr(\ba)=\Tr(\bb)=1+\delta m$, hence by Lemma~\ref{perm}, part (i), we have $s,t\in \bB_{1+\delta m}$
($s,t\neq 0,1,\infty$) such that
\[\ba=\phi(s)=\frac {1} {s+s^{-1}}\; {\rm and}\;\; \bb=\phi(t)=\frac {1}
{t+t^{-1}},\] where the map $\phi$ is defined before the statement
of Lemma~\ref{perm}. Plugging these into (\ref{1-to1equ5}) and
applying Lemma~\ref{zsumexp}, part (ii), we have
\[
s^{\sigma -1}+ s^{1-\sigma}+\theta(s+s^{-1})^{\sigma +1}=t^{\sigma
-1}+ t^{1-\sigma}+\theta(t+t^{-1})^{\sigma +1},
\]
that is,
\[\phi(w_{\theta}(s))=\phi(w_{\theta}(t)),\]
where $w_{\theta}$, $\theta\equiv 0$ or $1$ (mod 2), is defined
before the statement of Lemma~\ref{perm}.

By Lemma~\ref{perm}, part (i), since the map $\phi$ is two-to-one from $\bB_e$ to $\bT_e$ ($e=0$ or 1), we have
$w_{\theta}(s)=w_{\theta}(t)$ or $w_{\theta}(s)=w_{\theta}(t)^{-1}=w_{\theta}(t^{-1})$.
By (\ref{congr}), if $\theta=0$ and $1+\delta m\equiv 1$ (mod 2), then $k\equiv 1$ (mod 2); also, if $\theta=1$,
then $1+\delta m\equiv 0$ (mod 2) implies that $m\equiv 1$ (mod 2), and $1+\delta m\equiv 1$ (mod 2) implies that
$m\equiv 1+k$ (mod 2). So by Lemma~\ref{perm}, part (ii) and (iii), the map $z\mapsto w_{\theta}(z)$ is a
permutation of $\bB_{1+\delta m}$. Therefore we have either $s=t$ or $s=t^{-1}$, both lead to $\ba=\bb$, hence
$a=b$, therefore $x=y$. This completes the proof. \eop

\vspace{0.1in}

\noindent{\bf Remark 1.} In the above proof that $H_{\alpha, \gamma}$ is injective on $\bT_e$ for $e=0$ and
$e=1$, different proofs were given for the two cases. However, it is not difficult to adapt the above proof given
for the case $e=1$ so that it works for both cases $e=0$ and $e=1$. To this end, we first define the translation
maps $\tau_v$ for $v=0,1$ by $\tau_v(x)=x+v$. Now choose $\beta$ such that $k+\beta m\equiv 1$ (mod 2) (this is
possible since $k$ and $m$ are relatively prime), define $\delta$ as in Lemma~\ref{fgprop}, part (vi), and let
$\theta\in \{0,1\}$ satisfy $\theta \equiv \beta+\delta k$ (mod 2), so that (\ref{congr}) holds. Let $e \in
\{0,1\}$. We will in fact show that for all $z$ in $\bB_{e(1+\delta m)}$,
\begin{equation} \label{Hitt}
H_{\alpha,\gamma}(g_\beta(\tau_{\delta e}(\phi(z)))=
\tau_{\gamma e}(\phi(w_{\theta e}(z))).
\end{equation}
To see this, let $z\in \bB_{e(1+\delta m)}$. Now observe that $\by:=\phi(z)\in\bT_{e(1+\delta m)}$ by Lemma~\ref{perm}, part (i). So $y:=\by+\delta e=\tau_{\delta e}(\by)\in \bT_e$, and by Lemma~\ref{fgprop}, part (ii), we have $x:=g_{\beta}(y)\in\bT_e$. Furthermore, as a consequence of our choices for $\beta$, $\delta$, and $\theta$, we have by Lemma~\ref{perm}, part (vi), that $f_\alpha(g_\beta(y))=\by$, and by Lemma~\ref{perm},
part (vii) with $\lambda=\delta$, that $g_\beta(y)= \theta e + g_0(\by)$. Also, if $z\notin\{0,\infty\}$, we can
conclude from Lemma~\ref{zsumexp} that $\theta e + g_0^2(\phi(z))=\phi(z)^{\sigma+1} / \phi(w_{\theta e}(z))$,
and it is easily verified that this equation also holds when $z\in\{0,\infty\}$. Using these observations, we
conclude that \beqa
H_{\alpha,\gamma}(x) &=& H_{\alpha,\gamma}(g_{\beta}(y))\\
        &=& \gamma e+ f_\alpha(g_\beta(y))^{\sigma+1}/g^2_\beta(y)\\
        &=& \gamma e+ \by^{\sigma+1}/(\theta e + g^2_0(\by))\\
        &=& \gamma e + \phi(z)^{\sigma+1}/
                \biggl(\phi(z)^{\sigma+1} / \phi(w_{\theta e}(z)) \biggr)\\
        &=& \gamma e + \phi(w_{\theta e}(z)),
\eeqa
that is, (\ref{Hitt}) holds.

Now by Lemma~\ref{perm}, part (ii) and (iii) and by (\ref{congr}), the map $w_{\theta e}$ is a permutation on
$\bB_{e(1+\delta m)}$. Moreover, by Lemma~\ref{perm}, part (i), $\phi$ maps $\bB_{e(1+\delta m)}$ two-to-one onto
$\bT_{e(1+\delta m)} = \bT_{e(r+\alpha m)}$ (see Lemma~\ref{fgprop}, part (vi)), and this set is in turn mapped
one-to-one onto $\bT_{e(r+(\alpha+\gamma)m)}$ by the map $\tau_{\gamma e}$. So the composition map $z\mapsto
\tau_{\gamma e}(\phi(w_{\theta e}(z)))$ in the right-hand side of (\ref{Hitt}) is two-to-one from
$\bB_{e(1+\delta m)}$ onto $\bT_{e(r+(\alpha+\gamma)m)}$. On the other hand, $\phi$ maps $\bB_{e(1+\delta m)}$
two-to-one onto $\bT_{e(1+\delta m)}$, the map $\tau_{\delta e}$ maps this set one-to-one onto $\bT_e$, and
$g_\beta$ is a permutation on $\bT_e$, so the composition map $z\mapsto H_{\alpha,\gamma}(g_\beta(\tau_{\delta
e}(\phi(z))))$ is two-to-one if and only if $H_{\alpha,\gamma}$ is one-to-one on $\bT_e$. Combining these two
observations, we conclude that $H_{\alpha,\gamma}$ is one-to-one on $\bT_e$ for both $e=0$ and $e=1$.

\noindent{\bf Remark 2.} In the case where $\gamma=1$ and $m$ is odd, if $r+(\alpha+\gamma)m\equiv 1$ (mod 2),
then $r+\alpha m\equiv 0$ (mod 2), hence by Theorem~\ref{mainthm}, $H_{\alpha,0}(X)$ is not a PP of $\bF_q$. Yet,
by adding $\Tr(X)$ to $H_{\alpha,0}(X)$, we see that $H_{\alpha,1}(X)=\Tr(X)+H_{\alpha,0}(X)$ is a PP of $\bF_q$.

\noindent{\bf Remark 3.} When $k=1$ (so $\sigma=2$ and $r=1$), the map $H_{1,1}:\bF_q\rightarrow\bF_q$ fixes
$\bT_0$ elementwise and maps $x\in\bT_1$ to $x+1/x+1/x^2$. Therefore, by Theorem~\ref{mainthm} the map
$h:\bT_1\rightarrow \bT_1$ defined by $h(x)=x+1/x+1/x^2$ is a permutation of $\bT_1$. This fact was used in
\cite{kloos} to prove a Kloosterman sum identity.

\noindent{\bf Remark 4.} We give one more example to illustrate Theorem~\ref{mainthm}. Let $k, m$ be positive
integer such that $2k\equiv 1$ (mod $m$). Let $\sigma=2^k$. Then $\sigma^2\equiv 2$ (mod $2^m-1$). In this case,
we have $r=2$, and $$H_{0,0}(X)=X^{\sigma -1}+X^{2(\sigma -1)}+X^{\sigma^2-1}+X^{\sigma^2+\sigma-2},$$ and
$$H_{0,1}(X)={\rm Tr}(X)+X^{\sigma -1}+X^{2(\sigma -1)}+X^{\sigma^2-1}+X^{\sigma^2+\sigma-2}.$$ By
Theorem~\ref{mainthm}, $H_{0,0}$ maps $\bT_0$ bijectively to $\bT_0$, and $\bT_1$ bijectively to $\bT_0$; and
$H_{0,1}$ maps $\bT_0$ bijectively to $\bT_0$, and $\bT_1$ bijectively to $\bT_1$. In particular, $H_{0,1}(X)$, and hence also the polynomial 
$${\rm Tr}(X)+X^{\sigma -1}+X^{2(\sigma -1)}+X+X^{\sigma}$$ are PPs of $\bF_q$

\vspace{0.1in}

\noindent{\bf Acknowledgement:}  The research of the second author
was partially supported by NSA grant MDA 904-03-1-0095. The authors thank
an anonymous referee for his/her careful reading of the paper.

\end{document}